\documentclass[10pt,twoside,a4paper]{amsart}
\usepackage{amsmath}

\usepackage{amssymb}

\newcommand{\frow}{{}^\frown }

\newcommand{\dom}{\operatorname{dom}}
\newcommand{\ran}{\operatorname{ran}}
\newcommand{\depth}{\operatorname{depth}}
\newcommand{\height}{\operatorname{height}}

\renewcommand{\int}{\operatorname{int}}

\newcommand{\clop}{\operatorname{Clop}}
\newcommand{\Clop}{\operatorname{Clop}}

\newcommand{\Ult}{\operatorname{Ult}}

\newcommand{\Aut}{\operatorname{Aut}}
\newcommand{\aut}{\operatorname{Aut}}

\newcommand{\ext}{\operatorname{ext}}

\newcommand{\card}[1]{\mid\!\!{#1}\!\!\mid}

\newtheorem{thm}{Theorem}[section]
\newtheorem{corollary}[thm]{Corollary}
\newtheorem{lemma}[thm]{Lemma}

\theoremstyle{definition}
\newtheorem{defn}[thm]{Definition}
\newtheorem{claim}[thm]{Claim}

\begin{document}

\title[Some notes concerning homogeneity]{Some notes concerning the 
homogeneity of Boolean algebras and 
Boolean spaces}
\author{Stefan Geschke}
\author{Saharon Shelah}
\address[Geschke]{II.~Mathematisches Institut, Freie Universit\"at
Berlin, Arnimallee 3, 14195 Berlin, Germany.}
\email{geschke@math.fu-berlin.de}
\address[Shelah]{Institute of Mathematics, The Hebrew University of 
Je\-ru\-sa\-lem, 91904 Je\-ru\-sa\-lem, Israel and
Department of Mathematics, Rutgers University, New Brunswick, NJ 
08854, USA.}
\email{shelah@math.huji.ac.il}
\date{October 8, 2002}
\thanks{This article is [GeSh:811] in the second author's list of 
publications.  His research was supported by the German-Israeli 
Foundation for Scientific
Research \& Development and  the Edmund Landau Center for research
in Mathematical Analysis, supported by the Minerva Foundation
(Germany).}

\begin{abstract}
In this article we consider homogeneity properties of Boolean algebras 
that have nonprincipal ultrafilters which are countably generated. 

It is shown that a Boolean algebra $B$ is homogeneous if it is the union 
of countably generated nonprincipal ultrafilters and has a dense subset 
$D$ such that for every $a\in D$ the relative algebra $B\restriction 
a:=\{b\in B:b\leq a\}$ is isomorphic to $B$.  
In particular, the free product of countably many copies of an atomic 
Boolean algebra is homogeneous.  

Moreover, a Boolean algebra $B$ is homogeneous if it satisfies the 
following conditions:
\begin{itemize}\item[(i)] $B$ has a countably generated ultrafilter,
\item[(ii)] $B$ is not c.c.c., and 
\item[(iii)] for every $a\in B\setminus\{0\}$ there are 
finitely many automorphisms 
$h_1,\dots,h_n$ of $B$ such that $1=h_1(a)\cup\dots\cup h_n(a)$.
\end{itemize}

These results generalize theorems due to Motorov \cite{motorov} on the 
homogeneity of first countable Boolean spaces.  

Finally, we provide three constructions of first countable 
homogeneous Boole\-an spaces that are linearly ordered.  The first 
construction gives separable spaces 
of any prescribed weight in the interval $[\aleph_0,2^{\aleph_0}]$.
The second construction gives spaces of any prescribed weight in the 
interval $[\aleph_1,2^{\aleph_0}]$ that are not c.c.c. 
The third construction gives a space of weight $\aleph_1$ which is not 
c.c.c.~and which is not a continuous image of any of the previously 
described examples.  
\end{abstract}

\maketitle

\section{Introduction}
A topological space $X$ is {\em homogeneous} if for any two points $x,y\in 
X$ there is an autohomeomorphism of $X$ mapping $x$ to $y$.  Among the 
most obvious examples of homogeneous spaces are topological groups.  
In the case of topological groups, translations can be used to show their 
homogeneity.  

If we restrict our attention to zero-dimensional compact spaces, i.e., to 
Boolean spaces, topological groups are not interesting from the 
topological 
point of view since infinite compact zero-dimensional groups are all 
Cantor spaces, that is, they are homeomorphic to spaces of 
the form $2^\kappa$ where $\kappa$ is a cardinal (see 
\cite{hewitt} or \cite{hulanicki}).  
  
There is a surprising shortage of examples of homogeneous Boolean spaces 
other than Cantor spaces.  Interesting examples were provided 
by Maurice \cite{maurice}, who proved that for every indecomposable 
countable 
ordinal $\gamma$ the lexicographically ordered space $2^\gamma$ is 
homogeneous.  Here an ordinal $\gamma$ is {\em indecomposable} 
if $\gamma=\alpha+\beta$ with
$\beta>0$ implies $\beta=\gamma$.  If $\gamma>\omega$, then 
$2^\gamma$ ordered lexicographically does not satisfy the countable chain 
condition (c.c.c.) and therefore is not homeomorphic to a Cantor space.  

Cantor spaces and the lexicographically ordered spaces $2^\gamma$, 
$\gamma$ countable and indecomposable, have the property that not only 
the spaces themselves, but also their 
dual Boolean algebras are homogeneous.   
A Boolean algebra $B$ is {\em homogeneous} if for every $a\in 
B\setminus\{0\}$ 
the relative algebra $B\restriction a:=\{b\in B:b\leq a\}$ is isomorphic 
to $B$.
In general, there is no direct implication between the homogeneity of a 
Boolean algebra and the homogeneity of its Stone space. 
Van Douwen \cite{vandouwen} constructed a first countable homogeneous 
Boolean space whose 
dual Boolean algebra is not homogeneous.  And it is well known that the 
Boolean 
algebra $\mathcal P(\omega)/fin$ is homogeneous but its Stone space 
$\beta\omega\setminus\omega$ is not.  

However, the homogeneity of first countable Boolean spaces follows from 
the homogeneity of their dual Boolean algebra.  This was noticed 
independently by 
Motorov \cite{motorov} and Koppelberg \cite{koppelberg}.  
Motorov proved that the converse is also true in certain cases. 
He showed (in topological terms) that the homogeneity of a Boolean algebra 
follows from the homogeneity of its 
Stone space if the Boolean algebra is not c.c.c.~and every ultrafilter is 
countably generated.  Note that the last condition is equivalent to the 
first countability of the Stone space. 

The main tool in Motorov's argument is 
\begin{thm}\label{motorovmain} Let $B$ be a Boolean algebra such that 
every ultrafilter of 
$B$ is countably generated and $B$ has a dense subset $D$ such that for 
all $a\in D$, the algebra $B\restriction a$ is isomorphic to $B$.  Then 
$B$ is homogeneous.
\end{thm}  

Unfortunately, published proofs of Motorov's results seem to be 
unavailable.  We give the proofs of some generalizations of 
his theorems.  The main observation is that in Theorem \ref{motorovmain}
the assumption ``every ultrafilter of $B$ is countably generated'' can be 
weakened to ``$B$ is the union of countably generated ultrafilters'' 
(which is equivalent to the Stone space of $B$ having a dense set of 
points of countable character).
This easily implies that the free product of infinitely many copies of an 
atomic Boolean algebra is homogeneous.  Here a Boolean algebra $B$ is {\em 
atomic} if the atoms are dense in $B$, i.e., if the Stone space of $B$ 
has a dense set of isolated points.  

We also show that a Boolean 
algebra $B$ which is not c.c.c.~is homogeneous if it has at least one  
countably generated ultrafilter and the property that for all $a\in 
B\setminus\{0\}$ there are finitely many automorphisms $h_1,\dots,h_n$ 
of $B$ such that $1=h_1(a)\cup\dots\cup h_n(a)$.  The latter property is 
equivalent to the property that every point of the Stone space $X$ of 
$B$ has a dense orbit with respect to the natural group action of the 
group $\Aut(X)$ of autohomeomorphisms of $X$.  
 
Moreover, we  provide three  constructions 
leading to new examples of homogeneous 
Boolean spaces.  In all cases we obtain first countable spaces 
which are linearly ordered.  The first construction yields 
separable spaces of any prescribed weight in the interval 
$[\aleph_0,2^{\aleph_0}]$.  These spaces are constructed from nice 
suborders of $\mathbb R$.  Note that the space of countable weight is 
homeomorphic to $2^\omega$ since up to homeomorphism $2^\omega$ is the 
only Boolean space of countable weight without isolated points.  

The second 
construction uses an easy L\"owenheim-Skolem argument and 
gives homogeneous continuous images of the lexicographically 
ordered spaces $2^\gamma$, $\gamma<\omega_1$ indecomposable and 
uncountable.  The spaces obtained using this construction can have any 
prescribed weight in the interval 
$[\aleph_0,2^{\aleph_0}]$, and their cellularity equals their weight.
The third construction uses a linear order on an Aronszajn tree and yields 
a space of weight $\aleph_1$ which is not c.c.c.~and not a continuous 
image of any of the previously known examples of first countable 
homogeneous Boolean spaces.  

It should be noted that compact homogenous spaces which are linearly 
ordered have to be first countable (see \cite{bell} or \cite{kunen2}).  

\section{Generalizing Motorov's results}
As usual, the Stone space of a Boolean algebra $B$ is denoted by $\Ult(B)$ 
and the Boolean algebra of clopen subsets of a Boolean space $X$ is 
denoted by $\Clop(X)$.  In the following, we will frequently switch
between Boolean algebras and their Stone spaces, but our presentation will 
be mainly in topological terms.  

Note that a Boolean algebra $B$ is homogeneous if and only if every 
nonempty clopen subset of $\Ult(B)$ is homeomorphic to $\Ult(B)$.  

\begin{lemma}\label{ufhomofromhomo} Let $X$ be a Boolean space such that 
$\clop(X)$ is homogeneous.  If $x,y\in X$ are points of countable 
character, then there is an autohomeomorphism of $X$ mapping $x$ to 
$y$.  In particular, $X$ is homogeneous if it is first countable.  
\end{lemma}

\begin{proof}
Assuming that $X$ is infinite, it follows from the homogeneity of 
$\clop(X)$ that $X$ has no isolated points.  
Let $(A_n)_{n\in\omega}$ and $(B_n)_{n\in\omega}$ be clopen neighborhood 
bases of $x$ and $y$, respectively.  Since $x$ and $y$ are not isolated, 
we may assume that the sequences 
$(A_n)_{n\in\omega}$ and $(B_n)_{n\in\omega}$ are strictly decreasing.
We may also assume $A_0=B_0=X$.   
For each $n\in\omega$ let $C_n:=A_n\setminus A_{n+1}$ and 
$D_n:=B_n\setminus B_{n+1}$ and fix an homeomorphism $h_n:C_n\to D_n$.
It is easily checked that $h:=\{(x,y)\}\cup\bigcup_{n\in\omega}h_n$ is an 
autohomeomorphism of $X$ mapping $x$ to $y$.
\end{proof}

In order to apply Lemma \ref{ufhomofromhomo} we need a criterion for the 
homogeneity of Boolean algebras with first countable Stone spaces.  
A {\em $\pi$-base} of a topological space $X$ is a family $\mathcal F$ of 
open subsets of $X$ such that every nonempty open subset of $X$ includes a 
member of  $\mathcal F$.  A family of clopen subsets of a Boolean space 
$X$ is a $\pi$-base if and only if it is a dense subset of $\Clop(X)$.

\begin{lemma}\label{homofrompibase} 
Let $X$ be a Boolean space with a 
dense set of nonisolated points of 
countable character.  Then $\clop(X)$ is homogeneous if $X$ has a 
$\pi$-base consisting of clopen sets which are homeomorphic to $X$. 
\end{lemma}

\begin{proof}  We show that the nonempty clopen subsets of $X$ are 
pairwise homeomorphic.  Let $A$ be a nonempty clopen subset of $X$. 
Let $x\in A$ be a nonisolated point of countable character.  As in the 
proof of Lemma \ref{ufhomofromhomo}, there is a disjoint family 
$(A_n)_{n\in\omega}$ of nonempty clopen subsets of $A$ such that 
$A=\{x\}\cup\bigcup_{n\in\omega}A_n$.  

Inductively we define sequences $(C_n)_{n\in\{-1\}\cup\omega}$ and 
$(B_n)_{n\in\omega}$ as follows:  
Let $C_{-1}:=\emptyset$.   Let $n\in\omega$ and suppose we have already 
defined $C_{n-1}$.  Since the clopen subsets of $X$ which are homeomorphic 
to $X$ form a $\pi$-base of $X$, there is a clopen set $B_n\subseteq A_n$ 
such that $B_n$ is homeomorphic to $X\setminus C_{n-1}$.  With this 
choice,  $C_{n-1}\cup B_n\cong X$.  Let $C_n:= 
A_n\setminus B_n$.  

Now $$A\setminus\{x\}=\bigcup_{n\in\omega}(B_n\cup 
C_n)=B_0\cup\bigcup_{n\in\omega}(C_n\cup B_{n+1}).$$
By the choice of the $B_n$, $n\in\omega$, 
$\bigcup_{n\in\omega}(C_{n-1}\cup B_{n})$ is homeomorphic to the 
disjoint union of $\aleph_0$ copies of $X$.  
It follows that $A$ is the one-point compactification of the disjoint 
union of $\aleph_0$ copies of $X$.  Since $A$ was arbitrary, it 
follows that the nonempty clopen subsets of $X$ are pairwise homeomorphic. 
\end{proof}

Using Lemma \ref{homofrompibase} and Lemma \ref{ufhomofromhomo}, we can 
give an easy 
proof of the 
homogeneity of the lexicographically ordered spaces $2^\gamma$, 
$\gamma$ countable and indecomposable.  For every  
$\alpha<\gamma$ and every $x\in 2^{\alpha+1}$ the set 
$I_x:=\{y\in 2^\gamma:x\subseteq y\}$ is a clopen interval in $2^\gamma$.
By the indecomposability of $\gamma$, each $I_x$ is homeomorphic to 
$2^\gamma$.  Clearly, $$\{I_x:\alpha<\gamma\wedge x\in 
2^{\alpha+1}\}$$ is a $\pi$-base of $2^\gamma$.  Thus, $\Clop(2^\gamma)$ 
is homogeneous by Lemma \ref{homofrompibase}.  Now the homogeneity of 
$2^\gamma$ follows from Lemma \ref{ufhomofromhomo}. 

Another corollary of Lemma \ref{homofrompibase} gives information about 
free products of atomic Boolean algebras.  

\begin{corollary}\label{atomicprodhomo} Let $X$ be a Boolean space with 
a dense set of isolated points.  Then $\clop(X^\omega)$ is homogeneous.
\end{corollary}

\begin{proof} Let $D$ be the set of subsets of $X^\omega$ 
of the form $\{(x_0,\dots,x_{n-1})\}\times X^{\omega\setminus n}$ where 
each $x_i\in X$ is 
isolated.  Clearly, $D$ consists of clopen sets that are 
homeomorphic to $X^\omega$.  Since the isolated points are dense in $X$, 
$D$ is a $\pi$-base of 
$X^\omega$.  Those sequences $(x_i)_{i\in\omega}\in X^\omega$ for which 
each $x_i$ is isolated in $X$ form a dense subset of $X^\omega$, and each 
of these sequences is of countable character in $X^\omega$.  
Now it follows from Lemma 
\ref{homofrompibase} that $\clop(X^\omega)$ is homogeneous.
\end{proof}

Note that for every cardinal $\kappa$,  $\clop(X^\kappa)$ is isomorphic to 
the free product of $\kappa$ copies of $\clop(X)$.  
It is easily checked that $\clop(X^\kappa)$ is homogeneous if there is a 
cardinal $\lambda\leq\kappa$ such that $\clop(X^\lambda)$ is homogeneous.  
Therefore Corollary \ref{atomicprodhomo} implies that for a Boolean space 
$X$ with a dense set of isolated points, for every infinite cardinal 
$\kappa$ the Boolean algebra 
$\clop(X^\kappa)$ is homogeneous.
In  other words,  if $B$ is an atomic Boolean algebra, then every free 
product of infinitely many copies of $B$ is homogeneous.  

To proceed we need a technical lemma relating the cellularity of a compact 
space with many autohomeomorphisms to the cellularities of its nonempty open 
subsets.
For a topological space $X$ let $c(X)$ denote the cellularity of $X$.
Recall that $\aut(X)$ is the group of autohomeomorphisms of $X$.  For 
$x\in X$ the {\em $\aut(X)$-orbit} of $x$ is the set 
$\{h(x):h\in\aut(X)\}$. 

\begin{lemma}\label{cellularityopen} Let $X$ be compact and infinite.  If 
every $x\in X$ has a 
dense $\aut(X)$-orbit, then for every nonempty open subset $O$ of 
$X$ we have $c(O)=c(X)$. 
\end{lemma}

\begin{proof}
It is easily checked that all $\aut(X)$-orbits are dense in $X$ if and 
only if for every nonempty open set $O\subseteq X$, $\{h[O]:h\in\aut(X)\}$ 
covers $X$.  
Let $O\subseteq X$ be open and nonempty.  By the compactness of $X$, there 
are $n\in\omega$ and 
$h_1,\dots,h_n\in\aut(X)$ such that $X=h_1[O]\cup\dots\cup h_n[O]$.  

Let $\mathcal A$ be an infinite family of pairwise disjoint subsets of 
$X$.  For some $i\in\{1,\dots,n\}$, the set $\{A\in\mathcal A:A\cap 
h_i[O]\not=\emptyset\}$ is of size $\card{\mathcal A}$.  
It follows that $c(O)\geq\card{\mathcal A}$.  This implies $c(O)=c(X)$.  
\end{proof}

Now we have collected the necessary tools to show

\begin{thm}\label{nonccchomo}
Let $X$ be a Boolean space which is not c.c.c.~ and has a point of 
countable character. Suppose every $x\in X$ has a dense $\aut(X)$-orbit. 
Then $\clop(X)$ is homogeneous.   
\end{thm}  

\begin{proof}
Since $X$ is not c.c.c., $X$ is infinite.  Since every $\aut(X)$-orbit is 
dense in $X$, $X$ has no isolated points.  
Let $x_0\in X$ be a point of countable character.  Since the 
$\aut(X)$-orbit of $x_0$ is dense in $X$, $X$ has a dense set of  
points of countable character.  By Lemma \ref{homofrompibase}, it remains 
to show that $X$ has a $\pi$-base consisting of clopen sets which are 
homeomorphic to $X$.  

Let $(U_n)_{n\in\omega}$ be a neighborhood base of $x_0$ consisting of 
clopen sets.  
For every $n\in\omega$ there are $m\in\omega$ and 
$h_1,\dots,h_m\in\aut(X)$ such that $X=h_1[U_n]\cup\dots\cup 
h_m[U_n]$.  It follows that for each $n\in\omega$, $X$ is 
homeomorphic to a disjoint union of finitely many copies of clopen subsets
of $U_n$.  

Now let $O$ be a nonempty open subset of $X$.  By Lemma 
\ref{cellularityopen}, there is an uncountable family $\mathcal A$ of 
pairwise disjoint nonempty open subsets of $O$.
For every $A\in\mathcal A$ let $h_A\in\aut(X)$ be such that $h_A(x_0)\in 
A$.  $h_A$ exists since the $\aut(X)$-orbit of $x_0$ is dense.
For every $A\in\mathcal A$ there is $n(A)\in\omega$ such that 
$h[U_{n(A)}]\subseteq A$.  
Since $\mathcal A$ is uncountable, there is $n_0\in\omega$ such that for 
uncountably many $A\in\mathcal A$, $n(A)=n_0$. 
It follows that $O$ includes uncountably many pairwise disjoint open 
copies of  $U_{n_0}$.  But since $X$ is homeomorphic to a disjoint union 
of finitely many copies of clopen subsets of $U_{n_0}$, $O$ includes 
a clopen copy of $X$. 
\end{proof}

\begin{corollary}\label{noncccufhomo}  Let $X$ be a first countable 
Boolean 
space of 
uncountable cellularity.  If every point in $X$ has a dense $\aut(X)$-orbit,
then $\clop(X)$ and $X$ are both homogeneous.  In particular, $X$ is 
homogeneous if and only if $\clop(X)$ is.  
\end{corollary}  

\begin{proof} The homogeneity of $\Clop(X)$ follows immediately from 
Theorem \ref{nonccchomo}.  The homogeneity of $X$ now follows from Lemma 
\ref{ufhomofromhomo}.  
\end{proof}

\section{Examples of homogeneous Boolean spaces}
The homogeneous Boolean spaces we are going to construct will be Stone 
spaces of interval algebras of certain linear orders.  
As usual, if $(L,\leq)$ is a linear 
order,  we use $<$ to denote $\leq\setminus=$.  Similarly, if $<$ is 
transitive and irreflexive, we use $\leq$ to denote $<\cup=$.   

\begin{defn}  Let $(L,\leq)$ be a linear order.  The {\em interval 
algebra} 
$B(L)$ of $L$ is the subalgebra of $\mathcal P(L)$ generated by the 
intervals $[x,y)$, $x,y\in L$, $x<y$.     
\end{defn}  

Every element of $B(L)$
is a finite union of intervals of the form $[x,y)$, $x,y\in 
L\cup\{\infty\}$, $x<y$,  and of the form $(-\infty,x)$, $x\in 
L\cup\{\infty\}$ (see \cite{handbook}).  

The Stone space of an interval algebra $B(L)$ is homeomorphic to the 
linear order of initial segments of $L$ (see \cite{handbook}).  
Using this fact, we can characterize those linear orders whose interval 
algebras have a first countable Stone space.  

Call a subset $S$ of a linear order $L$ {\em coinitial} if for 
all $a\in L$ there is $b\in S$ such that $b\leq a$.  The {\em coinitiality} 
of $L$ is the least size of a coinitial subset of $L$, which is the same 
as the cofinality of the reversed order.

\begin{lemma}\label{firstcount} The Stone space of an interval algebra 
$B(L)$ is first 
countable if and only if every initial segment of $L$ has a countable 
cofinality and every final segment has a countable coinitiality.  
\end{lemma}

\begin{proof} Let $X$ be the set of initial segments of the linear order 
$L$.  $X$ itself is linearly ordered by $\subseteq$.  Suppose $X$ is first 
countable.  We show that every initial segment of $L$ is of 
countable cofinality.  The proof that every final segment of $L$ is of 
countable coinitiality is symmetric.  

By the first countability of $X$, 
for every $x\in X$, the set $\{y\in X:y\subsetneqq 
x\}$ is of countable cofinality.  Let $x\in X$ be nonempty and assume 
that $x$ has no last element.  Then the set 
$\{(-\infty,a]:a\in x\}$ 
is cofinal in $\{y\in X:y\subsetneqq x\}$.  Therefore, $\{(-\infty,a]:a\in 
x\}$ is of countable cofinality.  This implies that $x$ is of countable 
cofinality.  

Now suppose that every initial segment of $L$ is of countable cofinality 
and that every final segment of $L$ is of countable coinitiality.  
To show the first countability of $X$, it suffices to prove that for all 
$x\in X$ the following two conditions hold:  
\begin{enumerate}\item
If in $X$ there is no largest element below $x$, then $x$ is the first 
element of $X$ or there is a countable 
sequence in $X$ converging to $x$ from the left.  
\item If in $X$ there is no smallest element above $x$, then $x$ is the 
last element of $X$ or 
there is a countable sequence in $X$ converging to $x$ from the right.  
\end{enumerate}  
We show only the first condition since the proof of the second 
condition is 
symmetric.  Suppose that there is no largest element in $X$ 
which is below $x\in X$.  Assume further that $x$ is not the first 
element of $X$.  
Then $x$ as a subset of $L$ is nonempty and 
does not have a last element.  
By our assumption on $L$, there is a sequence $(a_n)_{n\in\omega}$ which 
is 
cofinal in $x$.   
We may assume that $(a_n)_{n\in\omega}$ is strictly increasing.  
The sequence $((-\infty,a_n])_{n\in\omega}$ of initial segments of $L$ 
converges to $x$ from the left.  
\end{proof}

Lemma \ref{firstcount} easily implies

\begin{corollary}\label{firstcor} If the linear order $L$ has no 
uncountable sequences (indexed by ordinals) which are strict\-ly 
increasing or strict\-ly 
decreasing, then the Stone space of $B(L)$ is first 
countable.  In particular, the Stone space of $B(L)$ is first countable if 
$L$ has a countable dense subset.     
\end{corollary}

\begin{proof}  
If $L$ has an initial segment of uncountable cofinality, then it has a 
strictly increasing sequence of length $\omega_1$.  Similarly, if $L$ has 
a final segment of uncountable coinitiality, then it has a strictly 
decreasing sequence of length $\omega_1$.  

If $L$ has a strictly increasing or strictly decreasing sequence of length 
$\omega_1$, then it is not c.c.c.~and therefore cannot have a countable 
dense subset.
\end{proof}

The following lemma provides an easy criterion for the homogeneity of 
an interval algebra.

\begin{lemma}\label{homogeneous}
Let $L$ be a linear order with the property that every two nonempty open 
intervals of $L$ are isomorphic.  Then $B(L)$ is homogeneous.  
\end{lemma}

\begin{proof} Since every two nonempty open intervals of $L$ 
are isomorphic, also the intervals of the form $[x,y)$ with $x\in L$, 
$y\in L\cup\{\infty\}$, $x<y$, and $(x,y)\not=\emptyset$ are pairwise 
isomorphic.  It follows that for all $n\in\omega$ and all 
$x_0,\dots,x_{2n+1}\in L\cup\{-\infty,\infty\}$ 
with $x_0<\dots<x_{2n+1}$ and $(x_0,x_{2n+1})\not=\emptyset$, 
$(x_0,x_1)\cup\bigcup_{i=1}^n[x_{2i},x_{2i+1})$ is 
isomorphic to $L$.  

Let $a\in B(L)\setminus\{\emptyset\}$.  Then for some $n\in\omega$ there 
are $x_0,\dots,x_{2n+1}\in L\cup\{-\infty,\infty\}$ 
with $x_0<\dots<x_{2n+1}$ such that either $x_0=-\infty$ and 
$a=(x_0,x_1)\cup\bigcup_{i=1}^n[x_{2i},x_{2i+1})$ or $x_0\in L$ and 
$a=[x_0,x_1)\cup\dots\cup[x_{2n},x_{2n+1})$.  
In either case, it is easily checked that $B(L)\restriction x$ is 
isomorphic to the interval algebra of 
$(x_0,x_1)\cup\bigcup_{i=1}^n[x_{2i},x_{2i+1})$.  
Since $(x_0,x_1)\cup\bigcup_{i=1}^n[x_{2i},x_{2i+1})$ is isomorphic to 
$L$, $B(L)\restriction a$ is isomorphic to $B(L)$.  
\end{proof}

Combining the information we have gathered so far we obtain

\begin{lemma}\label{crucial}  Let $L$ be a linear order such that every 
nonempty open interval of 
$L$ isomorphic to $L$.  Then $\Ult(B(L))$ is homogeneous if 
and only if $L$ has 
no strictly increasing or strictly decreasing sequences of length 
$\omega_1$.   In particular, $\Ult(B(L))$ is  
homogeneous if $L$ is separable.
\end{lemma}

\begin{proof}  
Since $L$ is isomorphic to every one of its nonempty open intervals,  
$B(L)$ is homogeneous by Lemma \ref{homogeneous}.  If $L$ has no 
strictly increasing or strictly decreasing sequences of uncountable 
length, then $\Ult(B(L))$ is first countable by Corollary \ref{firstcor} 
and homogeneous by Lemma \ref{ufhomofromhomo}.  Now suppose that 
$L$ has a strictly increasing or strictly decreasing sequence of 
uncountable length.  Then the Stone space of $B(L)$, being homeomorphic
to the linear order of initial segments of $L$, is not first countable and
therefore cannot be homogeneous, as mentioned in the introduction.
\end{proof}
  
It remains to construct linear orders with the properties required in 
Lemma \ref{crucial}.  We first construct some separable linear orders.

\begin{lemma}\label{geschkeorder}  For every cardinal 
$\kappa\in[\aleph_0,2^{\aleph_0}]$  
there is a separable dense linear order $L$ of size $\kappa$ without 
endpoints which 
is isomorphic to every one of its nonempty open intervals.
\end{lemma}

\begin{proof} 
We define two operations $f,f^-:\mathbb R^3\to\mathbb R$ as 
follows:  For $x,y\in\mathbb R$ with $x<y$ let $g_{x,y}$ be an order 
isomorphism from $\mathbb R$ onto $(x,y)$.  If $x,y,z\in\mathbb R$ are 
such that $x<y$, we let $f(x,y,z):=g_{x,y}(z)$.  Otherwise let 
$f(x,y,z):=0$.  If $x,y,z\in\mathbb R$ are such that $x<z<y$, we put 
$f^-(x,y,z):=g^{-1}_{x,y}(z)$.  Otherwise let $f^-(x,y,z):=0$.  

Now let $\kappa\in[\aleph_0,2^{\aleph_0}]$ be a cardinal.  Let $L$ be a 
subset of $\mathbb R$ of size $\kappa$ such that $\mathbb Q\subseteq L$ 
and $L$ is closed under $f$ and $f^-$.   
Then $L$ is separable since $\mathbb Q\subseteq L$.  
For $x,y\in L$ with $x<y$, $L$ is closed under the mappings 
$f(x,y,\cdot):z\mapsto f(x,y,z)$ and $f^-(x,y,\cdot):z\mapsto f^-(x,y,z)$.  
Therefore $f(x,y,\cdot)\restriction L$ is an isomorphism between $L$ and 
$(x,y)\cap L$. 
\end{proof}

\begin{corollary}  Let $\kappa\in[\aleph_0,2^{\aleph_0}]$ be a cardinal.
Then there is a homogeneous Boolean space of weight $\kappa$ which is 
separable 
and first countable. 
\end{corollary}  

\begin{proof} Let $L$ be a linear order of size $\kappa$ as in Lemma 
\ref{geschkeorder}.  By Corollary \ref{firstcor}, the Stone space of 
$B(L)$ is first countable.  Let $D\subseteq L$ be countable and dense.  
Then for all $a\in D$, the set of those $x\in B(L)$ which contain $a$ is 
an ultrafilter $F_a$.  Since $D$ is dense in $L$ and $B(L)$ does not 
contain any singletons by the density of $L$, $B(L)$ is the union 
of the $F_a$, $a\in D$, i.e., $B(L)$ is $\sigma$-centered.  This is 
equivalent to the separability of $\Ult(B(L))$.
Finally, $\Ult(B(L))$ is homogeneous by Lemma \ref{crucial}.
\end{proof}

In order to construct first countable homogeneous Boolean spaces of
a given weight in the interval $[\aleph_1,2^{\aleph_1}]$ that are not 
c.c.c., one can start with an indecomposable countable ordinal 
$\gamma>\omega$ and use the downward L\"owenheim-Skolem theorem 
to get a continuous image
of the lexicographically ordered space $2^\gamma$ with the right 
properties.

\begin{thm}  For every cardinal $\kappa\in[\aleph_0,2^{\aleph_0})$ 
there is a first countable homogeneous Boolean space $X$ of weight 
and cellularity $\kappa$.
\end{thm}

\begin{proof}
Let $\gamma>\omega$ be a countable indecomposable ordinal.  
Then the lexicographically ordered space $2^\gamma$ is a first countable
homogeneous Boolean space.  Since $2^\gamma$ is linearly ordered, it is 
the Stone space of an interval algebra (see \cite{handbook}).  
Let $L$ be a linear order with $2^\gamma\cong\Ult(B(L))$.

Since $\Ult(B(L))$ is first countable, in $L$ there is no strictly 
increasing or strictly decreasing sequence of length $\omega_1$.
The cellularity of 
$2^\gamma$ is $2^{\aleph_0}$, as is the weight. 
By Corollary \ref{noncccufhomo}, $B(L)$ is a homogeneous Boolean algebra.  

Let $\lambda$ be a sufficiently large cardinal and consider the structure
$(H_\lambda,\in)$ where $H_\lambda$ is the family of sets whose transitive 
closure is of size $<\lambda$.  Fix an antichain $\mathcal A\subseteq 
B(L)$ of size $2^{\aleph_0}$ and let $M$ be an elementary submodel of 
$(H_\lambda,\in)$ of size $\kappa$ such that $L,\kappa,\mathcal A\in M$ 
and $\kappa\subseteq M$.  

Let $B:= B(L)\cap M$.  By elementarity, $B=B(L\cap M)$. 
$L\cap M$ is a linear order without strictly increasing or strictly 
decreasing sequences of length $\omega_1$.  
Thus, $X:=\Ult(B)$ is first countable.
Again by elementarity, $B$ is homogeneous.  By Lemma 
\ref{ufhomofromhomo}, $X$ is homogeneous, too.  

Since $\kappa\subseteq M$, $B$ is of size $\kappa$, and so is $\mathcal 
A\cap M$.  It follows that the cellularity of $X$ is $\kappa$. 
\end{proof}

There is another interesting example of a first countable 
homogeneous Boolean space.  This one is constructed from a linear order
on an Aronszajn tree and is not the continuous image of any of the 
first countable homogeneous Boolean spaces
mentioned so far.  

Recall that a tree is {\em Aronszajn} if it is of height $\omega_1$, has 
only countable levels, and does not include an uncountable chain.  
If $T$ is a tree ordered by $\supseteq$, we will always assume that 
incomparable elements of $T$ are disjoint.

\begin{lemma}\label{shelahorder}  There is a dense linear order $L$ 
of size $\aleph_1$ without endpoints and with the 
following properties:
\begin{itemize}\item[(i)] $L$ has no strictly increasing or strictly 
decreasing sequences of length $\omega_1$.
\item[(ii)] $L$ is isomorphic to every one of its nonempty 
open intervals.
\item[(iii)] $L$ is not c.c.c.
\item[(iv)] $B(L)$ has a subset which is an Aronszajn tree (ordered by 
$\supseteq$).
\end{itemize}  
\end{lemma}

\begin{proof} 
For two functions $f$ and $g$ with the same domain we write $f=^*g$ if $f$ 
and $g$ agree on all but finitely many points of their common domain.  
Using the construction of an Aronszajn tree given in \cite{kunen}, we 
obtain a 
sequence $(f_\alpha)_{\alpha<\omega_1}$ such that 
each $f_\alpha$ is a 1-1 function from $\alpha$ into 
$\mathbb Q\cap(0,1)$ and 
for all $\alpha,\beta<\omega_1$ with $\alpha<\beta$,
$f_\alpha=^*f_\beta\restriction\alpha$.  

Now for each $\alpha<\omega_1$ let 
$$S_\alpha:=\{f\in{}^\alpha\mathbb 
Q:f\mbox{ is 1-1 and }f=^*f_\alpha\}$$ and 
$T_\alpha:=\bigcup_{\beta<\alpha}S_\beta$.
$T:=\bigcup_{\alpha<\omega_1}T_\alpha$ ordered by inclusion is an 
Aronszajn tree. 

We define a linear order on $T$.  Let $x,y\in T$ be such that $x\not=y$. 
If $x$ and $y$ are incomparable with respect to 
$\subseteq$, put $\Delta(x,y):=\min\{\nu\in\dom(x):x(\nu)\not=y(\nu)\}$ 
and let $x<y$ if $x(\Delta(x,y))<y(\Delta(x,y))$.
If $x\subseteq y$ and $\dom(x)=\alpha$, let $x<y$ if $y(\alpha)>\pi$ 
and $y<x$ if $y(\alpha)<\pi$. 

In other words, $T$ is ordered lexicographically after identifying each 
$x\in T$ with the function $x\frow\pi$ where $\frow$ 
denotes the concatenation of sequences and $\pi$ is identified with  
the sequence of length one with value $\pi\in\mathbb R$.

\begin{claim} $(T,\leq)$ is not c.c.c.  
\end{claim}

It is wellknown that $T$ ordered by reverse inclusion is not c.c.c.   
For example, for every $\alpha<\omega_1$ let $x\in S_{\alpha+1}$ be such that 
$x(\alpha)=0$.  
Since the $x_\alpha$, $\alpha<\omega_1$, are 1-1,  
$(x_\alpha)_{\alpha<\omega_1}$ is an antichain in 
$(T,\supseteq)$.  Note that for all $x\in T$, $\mathbb Q\setminus\ran(x)$ 
is infinite since otherwise the construction of the $f_\alpha$, 
$\alpha<\omega_1$, would break down at some point.  Thus, for every 
$\alpha<\omega_1$ 
there are $p_\alpha,q_\alpha\in\mathbb Q\setminus\ran(x_\alpha)$ with 
$p_\alpha<q_\alpha$.  Now for each $\alpha<\omega_1$, $x_\alpha\frow 
q_\alpha$ and $x_\alpha\frow p_\alpha$ are elements of $T$.  
Clearly, $\{(x_\alpha\frow p_\alpha,x_\alpha\frow 
q_\alpha):\alpha<\omega_1\}$ is an uncountable family of pairwise 
disjoint nonempty open intervals of $(T,\leq)$.  

\begin{claim}\label{cl1} In $(T,\leq)$ there is no strictly increasing or 
strictly decreasing sequence of length $\omega_1$.
\end{claim}

We only show that there is no strictly increasing sequence of length 
$\omega_1$ since the proof for decreasing sequences is the same.  
Suppose $(x_\alpha)_{\alpha<\omega}$ is a sequence in $T$ that is 
strictly decreasing (with respect to $\leq$).  
We construct a strictly increasing subsequence 
$(y_\alpha)_{\alpha<\omega_1}$ such that

\begin{enumerate}\item for all $\alpha,\beta<\omega_1$ with 
$\alpha<\beta$, $\dom(y_\alpha)<\dom(y_\beta)$ and

\item for all $\alpha,\beta,\gamma<\omega_1$, if $\alpha<\beta<\gamma$ and 
$\delta=\dom(y_\alpha)$, then 
$y_\beta\restriction(\delta+1)=y_\gamma\restriction(\delta+1)$.
\end{enumerate}

In order to construct $(y_\alpha)_{\alpha<\omega_1}$ first note that 
by the countability of every 
$T_\alpha$, $\alpha<\omega_1$, we can thin out 
$(x_\alpha)_{\alpha<\omega_1}$ 
and assume $\dom(x_\alpha)<\dom(x_\beta)$ for all 
$\alpha,\beta<\omega_1$ with $\alpha<\beta$.  

Suppose we have found $y_\alpha\in\{x_\nu:\nu<\omega_1\}$ with 
$\dom(y_\alpha)=\delta$.  The set 
$C:=\{x_\nu\restriction\delta+1:y_\alpha<x_\nu\wedge\nu<\omega_1\}$
is a countable linear order (ordered by $\leq\restriction C$).  
For every $x\in C$ let 
$$\ext(x):=\{x_\nu:x_\nu\restriction\delta+1=x\wedge\nu<\omega_1\}.$$ 
Clearly, $(\ext(x))_{x\in C}$ is a partition of 
$\{x_\nu:y_\alpha<x_\nu\wedge\nu<\omega_1\}$ into countably many convex 
sets.  Since for all $x,y\in C$ with $x<y$ all elements of $\ext(x)$ are 
below all elements of $\ext(y)$, for cofinality reasons $C$ has a last 
element $x$.  Now $\ext(x)$ is a final segment of 
$\{x_\nu:y_\alpha<x_\nu\wedge\nu<\omega_1\}$.  
Choose $y_{\alpha+1}\in\ext(x)$.  This finishes the successor step of the 
construction. 

Now suppose that 
$\alpha<\omega_1$ is a limit ordinal and $y_\beta$ has been 
chosen for all $\beta<\alpha$.  Since $\{x_\nu:\nu<\omega_1\}$ is of 
cofinality $\aleph_1$, the set $\{y_\beta:\beta<\alpha\}$ is bounded
in $\{x_\nu:\nu<\omega_1\}$.  Pick $y_\alpha\in\{x_\nu:\nu<\omega_1\}$
such that $y_\beta<y_\alpha$ for all $\beta<\alpha$.  It is easily checked 
that this construction yields a sequence with the desired properties.

To finish the proof of Claim \ref{cl1}, for each $\alpha<\omega_1$ let 
$\delta_\alpha:=\dom(y_\alpha)$ and consider the sequence 
$(y_{\alpha+1}\restriction\delta_\alpha)_{\alpha<\omega_1}$.  
By the choice of the $y_\alpha$, $\alpha<\omega_1$, for all 
$\alpha,\beta<\omega$ with $\alpha<\beta$ we have
$y_{\alpha+1}\restriction\delta_\alpha\subsetneqq
y_{\beta+1}\restriction\delta_\beta$, contradicting the fact that $T$ has 
no uncountable chains with respect to $\subseteq$.  

In order to prove that $(T,\leq)$ is isomorphic to each of its 
nonempty open intervals it suffices to show

\begin{claim}
For every $x\in T$, $(T,\leq)$ is isomorphic to $(x,\infty)$ and to
$(-\infty,x)$.  
\end{claim}

We only show $T\cong(x,\infty)$ since the proof of
$T\cong(-\infty,x)$ is symmetric.  Let $\delta:=\dom(x)+1$.

For each $y\in S_{\delta}$ let $\ext(y):=\{z\in 
T:\delta\subseteq\dom(z)\wedge z\restriction\delta=y\}$. As before, 
$\ext(y)$ is a convex subset of $T$.  For $y,z\in S_\delta$, 
$(\ext(y),\leq\restriction\ext(y))$ and 
$(\ext(z),\leq\restriction\ext(z))$ are isomorphic by the isomorphism 
mapping every $y^\prime\in\ext(y)$ to 
$z\cup y^\prime\restriction(\dom(y^\prime)\setminus\delta)$. 

As suborders of $(T,\leq)$,
$T_{\delta+1}$ and $T_{\delta+1}\cap(x,\infty)$ both are countable dense 
linear orders without endpoints.  $S_\delta$ is a dense and co-dense
subset of $T_{\delta+1}$ and $S_\delta\cap(x,\infty)$ is a dense and 
co-dense subset of 
$T_{\delta+1}\cap(x,\infty)$.  By the usual back-and-forth argument, 
there is an isomorphism $\varphi$ between $T_{\delta+1}$ and 
$T_{\delta+1}\cap(x,\infty)$ mapping $S_\delta$ onto 
$S_\delta\cap(x,\infty)$.  

For every $y\in S_\delta$ let $\varphi_y$ be an isomorphism between 
$\ext(y)$ and $\ext(\varphi(y))$.  Now we can construct an isomorphism 
$\psi$ between $T$ and $(x,\infty)$ by letting $\psi(y):=\varphi(y)$ for 
every $y\in T_\delta$ and $\psi(y):=\varphi_{y\restriction\delta}(y)$ for 
every $y\in T\setminus T_\delta$.  

Finally, we have to find an Aronszajn tree inside $B(T,\leq)$.
Let $T^\prime$ be the subtree of $T$ consisting 
of those elements of $T$ whose ranges are subsets of $(0,1)$.
As $(T,\subseteq)$,  $(T^\prime,\subseteq)$ is an Aronszajn tree.  
This is the place where we take advantage of the fact that the ranges of 
the $f_\alpha$, $\alpha<\omega_1$, are subsets of $(0,1)$.  
Mapping every $x\in T^\prime$ to the interval 
$[x\frow 0,x\frow 1)$ of $(T,\leq)$, we obtain an embedding of 
$(T^\prime,\subseteq)$ into $(B(T,\leq),\supseteq)$.  Note that this 
embedding maps incomparable elements of $T^\prime$ to disjoint members of 
$B(T,\leq)$.
\end{proof}

\begin{thm} There is a homogeneous Boolean space of weight
$\aleph_1$ which is first countable, not c.c.c., and not a continuous 
image of any of the lexicographically ordered spaces $2^\gamma$, 
$\gamma<\omega_1$.  
\end{thm}

\begin{proof}  Let $L$ be a linear order as in Lemma \ref{shelahorder}. 
Then $\Ult(B(L))$ is first countable by Corollary 
\ref{firstcor}.  $\Ult(B(L))$ is homogeneous by Lemma 
\ref{crucial} and $B(L)$ is of size $\aleph_1$ since $L$ is.   
$B(L)$ is not c.c.c.~since $L$ is not c.c.c.

Now suppose that $\Ult(B(L))$ is a continuous image of the 
lexicographically ordered space $2^\gamma$ for some $\gamma<\omega_1$.
Then $B(L)$ embeds into $\clop(2^\gamma)$.  By condition (iv) of Lemma 
\ref{shelahorder}, this implies that $\clop(2^\gamma)$ has a subset $T$ 
such that $(T,\supseteq)$ is an Aronszajn tree.  This contradicts

\begin{claim} Let $T\subseteq\Clop(2^\gamma)$ be such that $(T,\supseteq)$ 
is a tree whose levels are all countable.  Then $T$ is countable.
\end{claim}

Every element $a$ of $\clop(2^\gamma)$ can be uniquely written as a finite 
union of clopen intervals that are maximal convex subsets of $a$.
We may assume that $\gamma$ is infinite.  
If $a\in\clop(2^\gamma)$ is nonempty, let $\depth(a)$ be the 
least ordinal $\alpha<\gamma$ such that 
there are $x,y\in a$ with $x<y$ and $\Delta(x,y)=\alpha$ such that the 
closed interval $[x,y]$ 
is a maximal convex subset of $a$.
  
For every $a\in T$ let $\height(a)$ be the ordertype 
of $(\{b\in T:b\supsetneqq a\},\supseteq)$. 
We show that for every $a\in T$ and every $\alpha<\omega_1$ the 
following statement holds:

\begin{itemize}\item[$(*)_{a,\alpha}$] The set
$\{b\in T:b\subseteq a\wedge\depth(b)\leq\depth(a)+\alpha\}$ 
is countable.  
\end{itemize}

The claim follows from this since there 
is no $b\in T$ with $\depth(b)>\gamma$.  

We show $(*)_{a,\alpha}$ by induction on $\alpha<\omega_1$ simultaneously 
for all $a\in T$.  Let $a\in T$.  We start with proving $(*)_{a,1}$ since 
$(*)_{a,0}$ is trivial, i.e., there is no $b\in T$ with $b\subseteq a$ 
and $\depth(b)<\depth(a)$.  

If $b\in T$ is such that $b\subseteq a$ and $\depth(b)=\depth(a)$,  then 
there are $x,y\in b$ such that $x<y$, $[x,y]$ is a maximal convex subset 
of $b$, and $\Delta(x,y)=\depth(a)$.  It is easily 
checked that there is $z\in a$ such that either 
\begin{itemize}\item[a)] $z<x$, $[z,y]$ is a maximal convex subset of 
$a$, and $\depth(a)=\Delta(z,y)$ or
\item[b)] $y<z$, $[x,z]$ is a maximal convex subset of $a$, and 
$\depth(a)=\Delta(x,z)$.
\end{itemize}
Suppose that $\{b\in T:b\subseteq 
a\wedge\depth(b)=\depth(a)\}$ is uncountable.  Then there is $p\in a$ such 
that for uncountably many $b\in T$ with $b\subseteq a$, $p$ occurs as $y$ 
in a) or as $x$ in b).   This implies that $\{b\in T:b\subseteq a\wedge 
p\in b\}$ is uncountable.  By our assumption on trees ordered by 
$\supseteq$, any two elements of 
$T$ are either disjoint or comparable, and thus $\{b\in T:b\subseteq 
a\wedge p\in b\}$ is a chain in $T$.  But this contradicts the fact that 
$\Clop(2^\gamma)$ does not include any uncountable wellordered 
chain.  This 
finishes the proof of $(*)_{a,1}$.  

Now let $\alpha=\beta+1$ for some $\beta<\omega_1$ and suppose we have 
$(*)_{b,\beta}$ for all $b\in T$.  Let $a\in T$.  By $(*)_{a,\beta}$, 
there are only countably many $b\in T$ with $a\supseteq b$ and 
$\depth(b)<\depth(a)+\beta$.   Let $\delta<\omega_1$ be a bound for the 
heights of such $b$.  If $b\in T$ is minimal (with respect to the order 
$\supseteq$ on $T$) with $a\supseteq b$ and  
$\depth(b)=\depth(a)+\beta$, then $\height(b)\leq\delta+1$.  
It follows that there are only countably 
many $b\in T$ that are minimal with  $a\supseteq b$ 
and $\depth(b)=\depth(a)+\beta$.  Since $(T,\supseteq)$ is a tree,
every $b\in T$ with $a\supseteq b$ and $\depth(b)=\depth(a)+\beta$ is 
above a minimal element of $T$ with these properties.  
Applying $(*)_{b,1}$ for every minimal $b\in T$ with
$a\supseteq b$ and $\depth(b)=\depth(a)+\beta$, we obtain 
$(*)_{a,\beta+1}$.  

Finally suppose that $\alpha<\omega_1$ is a limit ordinal and for all 
$a\in T$ and all $\beta<\alpha$ we have $(*)_{a,\beta}$.  Then for all 
$a\in T$ the set
\begin{multline*}\{b\in T:b\subseteq 
a\wedge\depth(b)<\depth(a)+\alpha\}=\\\bigcup_{\beta<\alpha}\{b\in 
T:b\subseteq\alpha\wedge\depth(b)<\depth(a)+\beta\}\end{multline*}
is countable, which 
shows $(*)_{a,\alpha}$.  
\end{proof}

\end{document}